\documentclass[9pt,notitlepage,oneside]{amsart}

\usepackage{amssymb}
\usepackage{amsthm}
\usepackage{booktabs}
\usepackage{float}
\usepackage{hyperref}

\theoremstyle{plain}

\theoremstyle{definition}

\def\eq{\quad\Longrightarrow\quad}
\def\eqq{\quad\Longleftrightarrow\quad}
\newcommand{\disp}[1]{\displaystyle{{#1}}}

\newcommand{\nprime}{\,\text{prime}\,}

\newcommand{\field}[1]{\mathbb{#1}}

\begin{document}

\title{Prove or Disprove \\ 100 Conjectures from the OEIS}
\date{\today}
\author{Ralf Stephan}
\begin{abstract}
Presented here are over one hundred conjectures ranging from easy to 
difficult, from many mathematical fields. 
I also briefly summarize methods and tools that have led to this 
collection.\end{abstract}

\maketitle
\parindent0pt
\parskip5pt

\begin{flushright}
{\sl
Dedicated to all contributors to the OEIS, \\
on occasion of its 100,000th entry.}
\end{flushright}

The On-Line Encyclopedia of Integer Sequences \cite{b-o} is a database
containing the start terms of nearly 100,000 sequences (as of Autumn 2004),
together with formul\ae{} and references. It was reviewed by Sloane in 2003
\cite{b-s}. 
In the first section of this work, I describe how such a database can assist
with finding of conjectures, and the tools necessary to this end; section~2
then lists over one hundred so found propositions from many fields that 
have been checked numerically to some degree but await proof. Finally, I
conclude in section~3, giving a webpage that follows the status of the assertions. An appendix provides links to corresponding OEIS entries.

\section{Finding connections}
With a string of numbers representing an integer sequence, a conjecture appears 
already whenever a possible formula for a prefix of the sequence is found, or when a 
transformation is discovered that maps from one prefix to that of another sequence.
The OEIS consists just of such prefixes.
While the online interface to the OEIS allows searching for simple patterns
in the database, its `superseeker' eMail service is able to find formul\ae{}
for several types of sequences directly, by computation, or indirectly, by
applying transformations and comparing with the entries. However, this way,
it is not possible to do bulk searches of the whole or parts of the database.

To help with systematic work, Neil Sloane, the creator/editor/maintainer of 
the OEIS,  offers a file where only the 
numbers are collected, and this file served as input to several
computer programs I wrote over two years of work as associate editor. 
The methods can be divided according to program usage and intention:

\begin{itemize}
\item using a C program, I extended the database's~$10^5$ 
sequences with their first and second differences, 
with subsequent lexicographical sort and visual inspection using a simple
scrolling program (like Unix~\texttt{less}). This method sounds awkward but is 
not very much so since close matches show clearly through all noise. It yielded
several conjectures, among many `trivial' identities and a few hundred
duplicates.
\item for frequent offline lookups, I built a further extended database
by applying many transformations (bisections, pairwise sums, part-gcd, odd 
part, to state the more non-obvious) to a core set of OEIS sequences. 
The resulting file with
370MB size was searched with Unix \texttt{grep} when needed and using several
different strategies, and this yielded a good part of the conjectures of section~2, not to speak of another lot of identities that were included for
reference in the database.
\item special scans of the OEIS numbers, like that of Plouffe\cite{b-pl}, searching for sequences of specific 
kind, most notably for C--finite sequences, using my own implementation of 
\texttt{guessgf}\cite{b-bp} in Pari\cite{b-pari}, 
and for \emph{bifurcative}\footnote{formerly called
divide-and-conquer.} sequences. It has to be noted that scanning for specific
types of sequences is nontrivial and cannot be fully automatized---even in
the case of C--finiteness---, one still has to check for false positives.
From the scans, I have only included those conjectures that seemed the
most surprising to me; everything else served as immediate improvement of 
the information contained in the OEIS.
\end{itemize}
Comparing the quantities of output from the methods, simple transformations 
resulted in the most hits, I also had the impression that C--finiteness is quite
common: about~12 per cent of sequences have the property. On the other hand, to squeeze out more identities from
the database, ever more specialized transformations/scans would be necessary, 
with less and less results \emph{per method}. There is definitely a kind of 
fractality to it, perhaps reflecting expertise structure of OEIS submitters.

\section{The conjectures}
From the very start, Neil Sloane
had encouraged all kind of integer sequences to be submitted. The decision has 
lead to a statistically representative snapshot of sequences from all 
mathematical fields  where integers occur. A similar mixture is visible 
in my following conjectures which deal with power series expansions, 
number theory, and enumerative combinatorics, as well as additive and 
combinatorial number theory and, due to my special interest, bifurcative 
and other nonlinear sequences. This is also roughly the order in which
the statements are presented. I refrained from cluttering the presentation
with the respective A--numbers, also to encourage independent 
recalculation.

Please note that, throughout the section,
I use a set of commonly known abbreviations that are listed in Table~1 and
some of whom were introduced by Graham et al.\cite{b-gkp}.

\floatstyle{boxed}
\restylefloat{table}
\begin{table}[t]
\begin{center}
\begin{tabular}{cl}
$p$ & a prime\\
$[\text{$n$ even}]$ & 1 if $n$ even, 0 otherwise\\
\# & number of, cardinality $\ldots$\\
$(a,b)$ & greatest common divisor of $a,b$\\
\texttt{lcm}$\,(a,b)$ & least common multiple of $a,b$\\
$\lfloor x\rfloor$ & floor, greatest integer $\le x$\\
$\lceil x\rceil$ & ceiling, smallest integer $\ge x$\\
$a\,|\,b$ & $a$ divides $b$\\
$\{n\,\big|\;A\}$ & the set of all $n$ with property $A$\\
$d(n)$ & $\tau(n),\sigma_0(n)$, number of divisors of~$n$\\
$\phi(n)$ & Euler totient/phi function\\
$\sigma(n)$ & $\sigma_1(n)$, sum of divisors of $n$\\
$F_n$ & $n$-th Fibonacci number, with $F_0=0$\\
$P(x)$ & a polynomial in $x$\\
$[x^n]\,P(x)$ & $n$-th coefficient of polynomial or power series~$P$\\
$\big[\frac{x^n}{n!}\big]\,f(x)$ & $n$-th coefficient of Taylor expansion of~$f$\\
$\lg x$ & base-2 logarithm of $x$\\
$v_2(n)$ & dyadic valuation of~$n$, exponent of 2 in $n$\\
$e_1(n)$ & number of ones in binary representation of $n$\\
Res & resultant\\
\end{tabular}
\end{center}
\caption{Symbols and abbreviations.}
\end{table}

\subsection{Easy start: special functions, binomials, and more}

\begin{equation}(1-4x)^{3/2}=1-6x+\sum_{n>1}\frac{12(2n-4)!}{n!(n-2)!}x^n.\end{equation}
\begin{equation}(1+x^2C(x)^2)C(x)^2=\sum_{n\ge0}\frac{6n(2n)!}{n!(n+1)!(n+2)}x^n, \quad C(x)=\frac{1-\sqrt{1-4x}}{2x}.\end{equation}
Please show how to get the \emph{reduced} numerator/denominator in:
\begin{equation}(1-x)^{1/4}=\sum_{n\ge0}\frac{\prod_{k=1}^n(5-4k)2^{v_2(k)}/k}{2^{3n-e_1(n)}}x^n.\end{equation}
\begin{equation} [x^{n}]\,P_{n+2}(x) = \frac{1}{2^{n+2}}(n+1)\binom{2n+2}{n+1},\quad\text{($P_n$ the Legendre polynomials).}\end{equation}
\begin{equation}\sum_{k=0}^n-v_2\big(\,[x^{2k}]\,P_{2n}(x)\big)=2n^2+2n-2\sum_{i=0}^ne_1(i).\end{equation}
\begin{equation}B(2n,\tfrac12) = \frac{a}{b},B(2n,\tfrac14) = \frac{c}{d} \eq a=c, \quad\text{$B(n,x)$\ Bernoulli polynomials}.\end{equation}
\begin{equation}\arctan(\tanh x\tan x) = \sum_{n\ge0}(-1)^n2^{6n+2}(2^{4n+2}-1)\frac{B_{4n+2}}{2n+1}\cdot\frac{x^{4n+2}}{(4n+2)!}.\end{equation}
\begin{equation}\tau(2^n)=[x^{2^n}]\,x\prod_{k\ge1}(1-x^k)^{24}=[x^n]\frac{1}{2048x^2+24x+1}.\end{equation}
\begin{equation}7\;|\;\Big[\frac{x^{6k+4}}{(6k+4)!}\Big]\,\frac{1}{2-\cosh(x)}.\end{equation}
\begin{equation}4\;|\;\Big[\frac{x^{6k+4}}{(6k+4)!}\Big]\,\exp(\cos x-1),\quad11\;|\;\Big[\frac{x^{10k}}{(10k)!}\Big]\,\exp(\cos x-1)\ldots\end{equation}
\begin{equation}\text{Res}\,(x^n-1,4x^2-1)=4^n-2^n-(-2)^n+(-1)^n.\end{equation}
\begin{equation}\left\{n\,\big|\;n^2-1\,|\,\tbinom{2n}{n}\right\} \setminus \{2\} \subset\left\{n\,\big|\;n!(n-1)!\,|\,2(2n-3)!\right\}.\end{equation}
\begin{equation}\sum_{k=0}^n(k+1)\sum_{l=0}^k2^l\binom{k}{l}\binom{n-k}{l}=[x^n]\,\frac{1-x}{(1-2x-x^2)^2}.\end{equation}
\begin{equation}\sum_{k=0}^{n+1}(k+1)\bigg[\binom{2n+1}{k}-\binom{2n+1}{k-1}\bigg]=\frac{n+2}{2}\binom{2n+2}{n+1}-4^n.\end{equation}
\begin{equation}a(n)=\sum_{k=0}^n\bigg[\binom{n}{k}\mod2\bigg]\cdot2^k \eq a_{2n+1}=3a_{2n}.\end{equation}
\begin{equation}\sum_{k=0}^{\lfloor n/2\rfloor}D^k\binom{n}{2k+1} = [x^n]\frac{x}{1-2x+(1-D)x^2}.\end{equation}
\begin{multline}\left(\binom{2n}{n},\binom{3n}{n},\dotsc,\binom{(n-1)n}{n}\right)=1 \eqq \\
n+1=\sum_ip_i^{e_i} \;\wedge\; m = \max_ip_i^{e_i} \;\wedge\; \frac{n+1}{m} > m.\end{multline}
\begin{equation}\{a_n\,\big|\;\text{Least term in period of cont.~frac. of $\sqrt{a_n}=20$}\}=100n^2+n.\end{equation}
Define PCF$(n)$ the period of the continued fraction expansion for $\sqrt{n}$.
Then
\begin{equation}\text{PCF}\,(n)=\,\text{PCF}\,(n+1)\equiv1\bmod2 \eq n\equiv1\bmod24.\end{equation}
\begin{equation}\text{The largest term in the periodic part of the cont.~frac. of $\sqrt{3^n+1}$ is $2\cdot\big\lfloor(\sqrt3)^n\big\rfloor$.}\end{equation}
The numerators of the continued fraction convergents to~$\sqrt{27}$ are
\begin{equation}[x^n]\,\frac{5+26x+5x^2-x^3}{1-52x^2+x^4}.\end{equation}

\subsection{Classical number theory}

\begin{equation}n=5^i11^j \eq n|\sum_{k=1}^{10}k^n.\end{equation}
\begin{equation}n+1\,|\,d(n!^n).\end{equation}
Group the natural numbers such that the product of the terms of the $n$-th 
group is divisible by~$n!$. Let~$a_n$ the first term of the~$n$-th group. Then
\begin{equation}a_n=\bigg\lfloor\frac{(n-1)^2}{2}+1\bigg\rfloor.\end{equation}
\begin{equation}\#\{\text{cubic residues mod $8^n$}\}=\frac{4\cdot8^n+3}{7}.\end{equation}
\begin{equation}\text{lcm}\,(3n+1,3n+2,3n+3)=\tfrac34(9n^3+18n^2+11n+2)(3+(-1)^n).\end{equation}
(26a) Let $a_n=\prod_{k=1}^n\text{lcm}\,(k,n-k+1)$. Then~$a_n=n^2(n-1)!^2$ for~$n$ 
even, $n+1$ prime. Also, if $n$ is odd and~$>3$, $2(n+1)a_n$ is a perfect 
square, the root of which has the factor~$\frac12n(n-1)((n-1)/2)!$. 
\begin{equation}\left\{\min(x)\,\big|\;p\,|\,px-x-1\right\}=p-1,\quad p=\nprime(n).\end{equation}
\begin{equation}\emptyset=\{n\,\big|\;\forall\,i,j,\,0<i<j<n,\,n>3:\;2^n+2^i+2^j+1 \ \text{is composite}\}.\end{equation}
\begin{equation}p\equiv5\bmod12\eqq\emptyset=\left\{x\,\big|\;x^4\equiv9\bmod p\right\}.\end{equation}
\begin{equation}d(n)=d(n+1)=\cdots=d(n+6) \eq n\equiv5\bmod16.\end{equation}
\begin{equation}\frac{\sigma(2n)}{\sigma(n)}=\frac{4\cdot2^{v_2(n)}-1}{2\cdot2^{v_2(n)}-1}.\end{equation}
\begin{equation}\forall k>0:\quad\emptyset\not=\left\{x\,\big|\;\phi(x)=2^kp\right\} \eqq 2p+1\;\text{prime.}\end{equation}
\begin{equation}\left\{\text{composite $n$}\,\big|\;\phi(n+12)=\phi(n)+12 \wedge \sigma(n+12)=\sigma(n)+12
\right\}\eq n\equiv64\bmod72.\end{equation}
\begin{equation}\left\{n\,\big|\;\sigma(d(n^3))=d(\sigma(n^2))\right\}\eq n\equiv1\bmod24.\end{equation}
\begin{equation}\left\{n\,\big|\;\sigma(n)=2u(n)\right\} \eq n\equiv108\bmod216,
\quad u(n)=\sum_{\substack{d|n\\(d,n/\!d)=1}}d.\end{equation}
\begin{equation}\left\{n\,\big|\;t(n)=t(t(n)-n)\right\} = 
\left\{n\,\big|\;n=5\cdot2^k\vee n=7\cdot2^k,k>0\right\},t(n)=|\phi(n)-n|.\end{equation} 
\begin{equation}\left\{n\,\big|\;\phi(n^2+1)=n\phi(n+1)\right\} = \{8\} \cup
\left\{n\,\big|\;n^2+1=\nprime \wedge n+1=\nprime\right\}.\end{equation}
\begin{equation}\left\{n\,\big|\;|n-2d(n)-2\phi(n)-2|=2\right\} =
\{2,72\}\cup\left\{16p\,\big|\;p>2\right\}.\end{equation}
\begin{equation}\left\{\text{Local maxima of $\sigma(n)$}\right\}\subset
\left\{m\,\big|\;m=\sigma(l)\,\wedge\,l=\text{local maximum of $d(n)$}\right\}.\end{equation}
\begin{equation}F_n\bmod9\ \text{has period $24$.}\end{equation}
\begin{equation}(2^p-1,F_p)>1 \eq (2^{kp}-1,F_{kp})>1.\end{equation}
\begin{equation}(2^p-1,F_p)>1 \eq 8p\,|\,(2^{p}-1,F_p)-1.\end{equation}
\begin{equation}(2^p-1,F_p)>1\;\wedge\; p\not\equiv1\bmod10 \eq \frac{(2^{p}-1,F_p)-1}{8p}>1.\end{equation}
\begin{equation}\#\{k\,\big|\;F_k\,|F_n\}=d(n)-[\text{$n$ even}].\end{equation}
\begin{equation}\text{$k$ squarefree $\wedge\;\field{Q}\,(\sqrt{-k})$ has class number $n$} \eq \max k\equiv 19\bmod 24.\end{equation}
Define $\varsigma(n)$ the smallest prime factor of~$n$.
Let $a_n$ the least number such that the number of numbers $k\le a_n$ with 
$k>\varsigma(k)^n$ exceeds the number of numbers with $k\le\varsigma(k)^n$.
Then \begin{equation}a_n = 3^n + 3\cdot2^n + 6.\end{equation}

\subsection{Additive/combinatorial number theory}
\begin{equation}\#\Big\{m\,\big|\;m=[q^i]\prod_{k=1}^n\sum_{j=0}^kq^j\Big\}=[x^n]\frac{x(1-x+2x^2-2x^3+x^4)}{(1+x^2)(1-x)^3}.\end{equation}
\begin{equation}\#\Big\{m\,\big|\;m=[q^i]\prod_{k=1}^n1+\sum_{j=0}^kq^{2j+1}\,\wedge\,n>6\Big\}=n^2-3.\end{equation}
Let $A_{n\ge0}$ the number of distinct entries in the~$n\times n$ 
multiplication table, i.e., the number of distinct products $ij$ 
with~$1\le i,j\le n$, then $A_n$ goes~$\{0,1,3,6,9,14,\dots\}$. 
Let further $B$ the set of composite numbers
$>9$ that are not equal to the product of their \emph{aliquot divisors} 
(divisors of~$n$ without~$n$). Then
\begin{equation}A_n-A_{n-1}=\frac{\text{$n+$ its smallest divisor $>1$}}{2},
\ \text{if $n\not\in B$}.\end{equation}
Define the sequence $\langle a_n\rangle$ as $a_1=2$, $a_2=7$, and $a_n$ the
smallest number which is uniquely~$a_j+a_k$, $j<k$. The sequence starts
$\{2,7,9,11,13,15,16,17,19,21,25\dots\}$. Then
\begin{equation}a_n-a_{n-1}\ \text{has period $26$}.\end{equation}
Let $A$ the set of numbers whose cubes can be partitioned into two nonzero
squares, and $B$ the set of numbers that are the sum of two nonzero squares.
Then
\begin{equation}A=B.\end{equation}

\subsubsection{Sum-free sequences}
The sequence with start values $\{a_1,a_2,\dots,a_s\}$ and
further values $a_{n>s}$ satisfying ``$a_n$ is the smallest number $>a_{n-1}$
not of the form $a_i+a_j+a_k$ for $1\le i<j<k\le n$''. 
Then
\begin{equation}\text{Start with $\{1,2,3\}$}:\quad a_{n+6}-a_{n+5}=a_{n+1}-a_n,\quad \text{for $n>6$}.\end{equation}
\begin{equation}\text{Start with $\{0,1,2,3\}$}:\quad a_n=1\;\vee\; a_n\equiv 2,3 \bmod 8.\end{equation}
\begin{equation}\text{Start with $\{1,2,4\}$}:\quad a_n=\frac{26n-125-11\cdot(-1)^n}{4},\quad \text{for $n>14$}.\end{equation}
\begin{equation}\text{Start with $\{1,3,4\}$}:\quad a_{n+8}-a_{n+7}=a_{n+1}-a_n,\quad \text{for $n>7$}.\end{equation}
\begin{equation}\text{Start with $\{0,1,3,4\}$}:\quad a_n=1\;\vee\; a_n\equiv 3,4 \bmod10.\end{equation}
\begin{equation}\text{Start with $\{1,3,5\}$}:\quad a_{n+5}-a_{n+4}=a_{n+1}-a_n,\quad \text{for $n>6$}.\end{equation}
Let $a_1=1,a_2=2$ and $a_n$ the smallest number not of form~$a_i$, 
$a_i+a_{n-1}$, or~$|a_i-a_{n-1}|$. Then
\begin{equation}a_n=3n+3-2\cdot2^{\lfloor\lg(n+2)\rfloor},\quad \text{for $n>2$}.\end{equation}
Let $a_n=n$ for $n<4$ and $a_n$ the least integer $>a_{n-1}$ not of 
form~$2a_i+a_j$, $1\le i<j<n$. Then
\begin{equation}a_n=4n-10-(n\bmod2),\quad \text{for $n>3$}.\end{equation}

\subsubsection{Progression-free sequences.}
The sequence with start values $\{a_1,a_2,\dots,a_s\}$ and
further values $a_{n>s}$ satisfying ``$a_n$ is the smallest number $>a_{n-1}$
that builds no three-term arithmetic progression with any $a_k$, $1\le k<n$'' 
is well--defined. This peculiar definition was chosen to allow for any
start values. We abbreviate such a sequence with start values~$a,b,\dots$ 
as~$A_3(a,b,\dots,n)$ and present several conjectures for them.
\begin{equation}A_3(1,3,n)=A_3(1,2,n)+[\text{$n$ is even}]=\frac{3-n}{2}
+2\lfloor n/2\rfloor+\frac12\sum_{k=1}^{n-1}3^{v_2(k)}.\end{equation}
\begin{equation}A_3(1,4,n)=A_3(1,3,n)+[\text{$n$ is even}]+
[\text{$\lceil n/2\rceil$ is even}].\end{equation}
\begin{equation}A_3(1,7,n)=b_n+\sum_{k=1}^{n-1}\frac{3^{v_2(n)}+1}{2},\ 
\text{where $b_n=\{1,6,5,6,2,6,5,7,\dots\}_{n\ge1}$ has period 8}.
\end{equation}
\begin{multline}A_3(1,10,n)=b_n+\sum_{k=1}^{n-1}\frac{3^{v_2(n)}+1}{2},\\
\text{where $b_n=\{1,9,8,9,5,10,10,10,\dots\}_{n\ge1}$ has period 8}.
\end{multline}
\begin{multline}A_3(1,19,n)=b_n+\sum_{k=1}^{n-1}\frac{3^{v_2(n)}+1}{2},\\
\text{where $b_n=\{1,18,17,18,14,18,17,19,5,18,17,18,14,19,19,19\dots\}_{n\ge1}$ has period 16}.
\end{multline}
\begin{multline}\text{In general,}\quad A_3(1,m,n)=b_n+\sum_{k=1}^{n-1}\frac{3^{v_2(n)}+1}{2},\quad m=1+3^k\;\vee\;1+2\cdot3^k, k\ge0, \\
\text{where $b_n$ has period $P\le2^{\lfloor\frac{k+3}{2}\rfloor}$}.
\end{multline}
\begin{equation}A_3(1,2,3,n+2)=1+2^{\lfloor\lg n\rfloor}+
\sum_{k=1}^n\frac{3^{v_2(n)}+1}{2}.\end{equation}
The numbers $n$ such that the $n$-th row of Pascal's triangle contains
an arithmetic progression are
\begin{equation}n=19\,\vee\, n=\tfrac18\big[2k^2+22k+37+(2k+3)(-1)^k\big],k>0.\end{equation}

\subsection{Enumerative combinatorics}\hfill\break
The number of $n\times n$ invertible binary matrices $A$ such that $A+I$ is invertible is
\begin{equation}2^{\frac{n(n-1)}{2}}a_n, \quad\text{with\ }
a_n=\langle a_0=1,a_n=(2^n-1)a_{n-1}+(-1)^n\rangle. \end{equation}
Consider flips between the $d$-dimensional tilings of the unary zonotope $Z(D,d)$. Here the codimension $D-d$ is equal to 3 and $d$ varies. Then the number of flips is
\begin{equation}a(d)=(d^2+11d+24)2^{d-1}. \end{equation}
The sequence $a_n$ shifts left twice under binomial transform and is described by both $\langle a_0=a_1=1,a_n=\sum_{k=0}^{n-2}\binom{n-2}{k}a_k\rangle$ 
and (formally)
\begin{equation}\sum_{n\ge0}a_nx^n=\sum_{k\ge0}\frac{x^{2k}}{(1-kx)(1-x-kx)\prod_{m=0}^{k-1}(1-mx)^2}. \end{equation}
The number of self-avoiding closed walks, starting and ending at the origin,
of length~$2n$ in the strip~$\{0,1,2\}\times\field{Z}$ is
\begin{equation}\frac{1}{125}\Big\{(315n-168)2^{n-2}+(-1)^{\lfloor n/2\rfloor}\big[55\lfloor n/2\rfloor+78-(135\lfloor n/2\rfloor+36)(-1)^n\big]\Big\}, n>1. \end{equation}
The number of non-palindromic reversible strings with $n$ beads of 4 colors is
\begin{equation}\begin{cases}\tfrac124^n-\tfrac122^n, &\text{$n$ even,}\\\tfrac124^n-2^n, &\text{$n$ odd.}\end{cases} \end{equation}
The number of non-palindromic reversible strings with $n-1$ beads of 2 colors (4 beads are black) is
\begin{equation}\begin{cases}\tfrac{1}{48}(n^4-10n^3+32n^3-38n+15), &\text{$n$ odd,}\\\tfrac{1}{48}(n^4-10n^3+32n^3-32n), &\text{$n$ even.}\end{cases} \end{equation}
The number of non-palindromic reversible strings with $n$ black beads and $n-1$ white beads is
\begin{equation}\begin{cases}\disp{\frac{1}{4}\left[\binom{2n}{n}-\binom{n}{n/2}\right]}, &\text{$n$ even,}\\\disp{\frac{1}{2}\left[\binom{2n}{n}-\binom{2n-1}{n-1}-\binom{n-1}{(n-1)/2}\right]}, &\text{$n$ odd.}\end{cases} \end{equation}
The number of necklaces of $n$ beads of 2 colors (6 of them black) is
\begin{equation}[x^n]\,\frac{x^6(1-x+x^2+4x^3+2x^4+3x^6+x^7+x^8)}{(1-x)^4(1+x)^2(1-x^3)(1-x^6)}. \end{equation}
The number of edges in the 9-partite Turan graph of order $n$ is
\begin{equation}[x^n]\,\frac{x}{(1-x)^2}\Big[\frac{1}{1-x}-\frac{1}{1-x^9}\Big]. \end{equation}
The number of binary strings of length $n$ that can be reduced to null by repeatedly removing an entire run of two or more consecutive identical digits is
\begin{equation}2^n-2nF_{n-2}-(-1)^n-1.\end{equation}
The number of nonempty subsets of $\{1,2,\dots,n\}$ in which exactly $1/2$ of the elements are $\le(n-1)/2$ is
\begin{equation}\binom{n}{\lfloor (n-1)/2\rfloor}-1.\end{equation}
The number of level permutations of $2n-1$ is
\begin{equation}\frac{(2n-1)!}{2^{2n-2}}\binom{2n-2}{n-1}.\end{equation}
The number of rooted trees with $n$ nodes and 3 leaves is
\begin{equation}\tfrac{1}{288}\big(6n^4-40n^3+108n^2-120n-41+9(-1)^n+32[(n\bmod3)+(n+1\bmod3)]\big).\end{equation}
Let $a_n$ the number of $(2\times n)$ binary arrays with a path of adjacent 
1's from the upper left corner to anywhere in right hand column. Then
\begin{equation}a_{n+2}=2P_n+5P_{n+1},\quad\text{$P_n=$ Pell numbers.}\end{equation}
The number of strings over $\field{Z}_3$ of length~$n$ with trace~0 and subtrace~1 is
\begin{equation}[x^n]\,\frac{x(-6x^4+6x^3)}{(1-3x)(1+3x^2)(1-3x+3x^2)}. \end{equation}
The number of strings of length $n$ over GF(4) with trace~0 and subtrace~0 is
\begin{equation}[x^n]\,\frac{x(-26x^3+13x^2-5x+1)}{(1-2x)(1-4x)(1+4x^2)}. \end{equation}
The number of symmetric ways to lace a shoe that has~$n$ pairs of eyelets, such that each eyelet has at least one direct connection to the opposite side, is
\begin{equation}\sum_{k=0}^nk!\binom{n}{k}F_{k+2}. \end{equation}
The number of minimax trees with $n$ nodes is~$2^n$ times the number of
labelled ordered partitions of a $2n$-set into odd parts, that is,
\begin{equation}2^n\,\bigg[\frac{x^{2n}}{(2n)!}\bigg]\,\frac{1}{1-\sinh x}.\end{equation}
The number of unlabeled alternating octopi with~$n$ black nodes and~$k$ white 
nodes has the g.f.
$$\sum_{k,n\ge1} \frac{\phi(k)}{k}\log\bigg(\frac{(1-x^ny^n)^2}{1-x^ny^n(3+x^n+y^n)}\bigg).$$
The conjecture is now that the number of those octopi with~$n$ black nodes and~$n$ white nodes (the diagonal of the above array) is
\begin{equation}-2+3\sum_{d|n}\frac{\phi(n/d)\binom{2d}{d}}{2n}.\end{equation}
Finally, the number of dimer tilings of the graph $S_k\times P_{2n}$ ($S_k$ the star graph on~$k$ nodes, $P_n$ the path with length $n$) is
\begin{equation}[x^n]\,\frac{1-x}{1-(k+1)x+x^2}.\end{equation}

\subsection{Nonlinear recurrences and other sequences}
\begin{equation}a_n=[x^n]\,\frac{-4x^5+x^4+x^3-3x^2-2x+6}{(1-x)(1-x-x^2-x^5)} \eqq 
\Big\langle a_0=6,a_1=10,a_n=\Big\lfloor\frac{a_{n-1}^2}{a_{n-2}}+
\frac12\Big\rfloor\Big\rangle.\end{equation}
\begin{equation}a_n=[x^n]\, \frac{-x^5+x^4-x^3+x^2-2x+3}{(1-x)(1-2x-x^3-x^5)} \eqq
\Big\langle a_0=3,a_1=7,a_n=\Big\lfloor\frac{a_{n-1}^2}{a_{n-2}}\Big\rfloor\Big\rangle.\end{equation}
\begin{equation}a_n=[x^n]\, \frac{-3x^5+2x^4+x^3-x^2-2x+4}{(1-x)(1-2x-x^2-2x^5)} \eqq
\Big\langle a_0=4,a_1=10,a_n=\Big\lfloor\frac{a_{n-1}^2}{a_{n-2}}\Big\rfloor\Big\rangle.\end{equation}
\begin{equation}a_n=[x^n]\, \frac{2x^3+x^2-4x+5}{-x^4+2x^2-3x+1} \eqq
\Big\langle a_0=5,a_1=11,a_n=\Big\lfloor\frac{a_{n-1}^2}{a_{n-2}}+
\frac12\Big\rfloor\Big\rangle.\end{equation}
\begin{equation}a_n=[x^n]\, \frac{3x^5+2x^4+x^3+4x^2-x+6}{-x^6-x^3+x^2-2x+1} \eqq
\Big\langle a_0=6,a_1=11,a_n=\Big\lfloor\frac{a_{n-1}^2}{a_{n-2}}+
\frac12\Big\rfloor\Big\rangle.\end{equation}
\begin{multline}\langle a_0=a_1=1, a_n=(|n-1-a_{n-1}|\bmod n-1)+(|n-1-a_{n-2}|\bmod n)\rangle \\ \eq a_{n+3}=a_n,\quad n>6.\end{multline}
\begin{equation}\langle a_0=x,a_{n+1}=a_n(a_n+1)\rangle \eq [x^{2^n-3}]\,a_k=\frac{2^{3k+2}-2^k}{3}.\end{equation}
\begin{equation}\bigg\langle \nu_0=\nu_1=1,\;\nu_n=\nu_{n-1}+3\nu_{n-2}\sum_{i=0}^{n-2}q^i\bigg\rangle \eq [q^1]\,\nu_n = [x^n]\frac{x^3(9x+3)}{(1-3x-3x^2)^2}.\end{equation}
\begin{equation}a_n=\sum_{k=1}^\infty\bigg\lfloor2\Big(\frac{\sqrt5+1}2\Big)^{n-k}\bigg\rfloor\eqq \sum_{n\ge0}a_nx^n=\frac{x(x^5+x^4-4x^2+3)}{(1-x)(1-x^2)(1-x-x^2)}.\end{equation}
\begin{equation}\bigg\langle a_0=1,a_{n+1}=\bigg\lfloor\frac{a_n}{\sqrt5-2}\bigg\rfloor\bigg\rangle \eqq \sum_{n\ge0}a_nx^n=\frac{1-x-x^2}{(1-x)(1-4x-x^2)}.\end{equation}
Let the sequence~$\langle a_n\rangle$ be defined such that $a_1=C$ 
and~$a_{n+1}=$ the smallest difference $>1$ between~$d$ and~$p/d$ for any 
divisor~$d$ of the partial product~$p=\prod_{k=1}^na_k$ of the sequence. Then
$a_n=\{19,18,29,27,9,\dots\}$ for $C=19$, and $a_n=\{21,4,5,13,8,2,\dots\}$ for
$C=21$ and
\begin{equation}a_n=3^k, C=19, n>3 \quad \text{and}\quad a_n=2^m, C=21, n>4.\end{equation}

\subsection{Binary representation, $k$-regular and bifurcative sequences}
\begin{equation}\#\big\{m\,\big|\;m=v_2\tbinom{n}{j},\,0\le j\le n\big\}=\lfloor\lg (n+1)\rfloor+1-v_2(n+1).\end{equation}
\begin{equation}e_1(m)\equiv0\bmod2 \eqq m\in\langle a_0=0,\;a_{2n}=a_n+2n,\;a_{2n+1}=-a_n+6n+3\rangle.\end{equation}
\begin{equation}e_1(m)\equiv1\bmod2 \eqq m\in\langle a_0=1,\;a_{2n}=a_n+2n,\;a_{2n+1}=-a_n+6n+3\rangle.\end{equation}
\begin{equation}a_n=n+[x^n]\frac{x}{1-x}\sum_{k\ge0}2^kx^{3\cdot2^k}\eqq \text{$a_n$ in binary does not begin $100$}.\end{equation}
\begin{multline}\langle a_0=0,a_1=1,a_{2n}=a_n,a_{2n+1}=a_{n+1}-a_n\rangle \;\wedge\; a_{3k}=0 \\\eqq \text{no adjacent 1s in binary of $k$}.\end{multline}
\begin{multline}\langle a_1=3,a_{2n}=4a_n-2n,a_{2n+1}=4a_n-2n+2^{\lfloor\lg(4n+2)\rfloor}\rangle \\ \eqq \text{binary of $a_n$ is binary of $n$ twice juxtaposed.}\end{multline}
\begin{equation}a_n=n\;\texttt{XOR}\;(n+m) \eqq a_n=[x^n]
\frac{P(x)}{(1-x)^2\prod_{k\ge0}1+x^{2^{e_k}}},\quad\sum_{k\ge0}2^{e_k}=m.
\end{equation}
\begin{multline}\langle a_0=a_1=0,a_{4n}=2a_{2n},a_{4n+2}=2a_{2n+1}+1,a_{4n+1}=2a_{2n}+1,a_{4n+3}=2a_{2n+1}\rangle \\ \eqq a_n=\,\{\text{Replace each pair of adjacent bits of $n$ by their mod 2 sum}\,\}.\end{multline}
\begin{equation}a_n=\sum_{k=1}^{n-1}k\,\texttt{AND}\,(n-k) \eqq 
\langle a_0=a_1=0,a_{2n}=2a_{n-1}+2a_n+n,a_{2n+1}=4a_n\rangle.\end{equation}
\begin{equation}a_n=\sum_{k=1}^{n-1}k\,\texttt{XOR}\,(n-k) \eqq 
\langle a_0=a_1=0,a_{2n}=2a_{n-1}+2a_n+4n-4,a_{2n+1}=4a_n+6n\rangle.\end{equation}
\begin{equation}a_n=\sum_{k=1}^{n-1}k\,\texttt{OR}\,(n-k) \eqq 
\langle a_0=a_1=0,a_{2n}=2a_{n-1}+2a_n+5n-4,a_{2n+1}=4a_n+6n\rangle.\end{equation}
\begin{multline}a_n=\#\{(i,j)\,|\,0\le i,j< n \,\wedge\,i\,\texttt{AND}\,j>0\} \\ \eqq \langle a_0=a_1=0,a_{2n}=3a_n+n^2,a_{2n+1}=a_n+2a_{n+1}+n^2-1\rangle.\end{multline}
\begin{multline}n=\sum_{k\ge0}2^ke_k\;\wedge\;a_n=\sum_{k\ge0}(-1)^ke_k\;\wedge\;|a_n|=3 \\ \eq n\in\left\{m\,\big|\;m=3k\,\wedge\,k=3i\,\wedge\,e_1(k)\equiv1\bmod2\right\}.\end{multline}
\begin{multline}n=\sum_{k\ge0}2^ke_k\;\wedge\;a_n=\sum_{k\ge0}(-1)^ke_k\;\wedge\;a_n=0 \\ \eq n=3m\;\wedge\;m\not\in\left\{k\,\big|\;k=3i\,\wedge\,e_1(k)\equiv1\bmod2\right\}.\end{multline}
\begin{equation}\langle a_0=0,a_{2n}=1-a_n,a_{2n+1}=-a_n\rangle\,\wedge\,a_{3k}=0 \eq
\text{$k$ in base-4 contains only $-1,0,1$.}\end{equation}
\begin{equation}\max\sum_{j=0}^n[x^j]\sum_{k\ge0}\frac{x^{2^k}}{1+x^{2^k}+x^{2^{k+1}}} = \lfloor\log_4n\rfloor+1.\end{equation}
Define the sequence $a_n$ by $a_1=1$ and $a_n=M_n+m_n$, where $M_n=
\max_{1\le i<n}(a_i+a_{n-i})$, and $m_n=\min_{1\le i<n}(a_i+a_{n-i})$. Let
further $b_n$ the number of partitions of~$2n$ into powers of~2 
(number of \emph{binary partitions}). Then
\begin{equation}m_n=\tfrac32b_{n-1}-1,\quad M_n=n+\sum_{k=1}^{n-1}m_n,
\quad a_n=M_{n+1}-1.\end{equation}
Let $a_n$ defined as the limit in the infinite of the sequence $b_n$, with $b_1=1$, $b_2=n$, and $b_{n+2}=\lceil\frac12(b_n+b_{n+1})\rceil$, and $c_n$ the number of ones in the base-$(-2)$-representation of~$n$. Then
\begin{equation}c_n=3a_{n+1}-2n-3\quad\text{and}\quad a_{n+1}-a_n=\langle d_{4n}=0,d_{2n+1}=1,d_{4n+2}=d_{n+1}\rangle.\end{equation}
Let $a_n$ the number of subwords of length~$n$ in the word generated by \texttt{a} $\mapsto$ \texttt{aab}, \texttt{b} $\mapsto$ \texttt{b}. Then
\begin{equation}\sum_{n\ge0}a_nx^n = 1+\frac{1}{1-x}+\frac{1}{(1-x)^2}\bigg(\frac{1}{1-x}-\sum_{k\ge1}x^{2^k+k-1}\bigg).\end{equation}

\vfill\pagebreak
\section{Conclusions}
Working over two years with the OEIS showed me that simple computer programs
suffice for many tasks; where I had to write programs myself, it was not
visible that the task could be fully automatized---mathematics is essentially
human. The work as editor was rewarding not only in itself but also in that
it yielded a huge collection of conjectures as byproduct. However, I expect
further gains in that regard as becoming ever more difficult as scans and
transformations have to become more specialized and complex.

I do not intend to work on proving the majority of propositions presented
here but I provide below a webpage giving the status of work done on them.
Given their number, it is quite possible that a few are already in the
literature. I hope the reader excuses my not researching these ones: it
is very difficult nowadays to access pay--only journals from outside
university.

\section{Acknowledgments}
I want to thank Elizabeth Wilmer, Robin Chapman, Jason Dyer, Ira Gessel, Mitch Harris,
Vladeta Jovovic, Nikolaus Meyberg, Luke Pebody, John Renze, and Lawrence Sze 
who pointed out several errors in the first versions of the file.


Status page: \href{http://www.ark.in-berlin.de/conj.txt}{\texttt{http://www.ark.in-berlin.de/conj.txt}}
\hfill\break
e-Mail: \href{mailto:ralf@ark.in-berlin.de}{\texttt{ralf@ark.in-berlin.de}}
\end{document}